\documentclass[a4paper,11pt,oneside,reqno,psamsfonts]{amsart}
\usepackage{a4wide}
\usepackage[T1]{fontenc} 
\usepackage[reqno]{amsmath}
\usepackage[psamsfonts]{amsfonts}
\usepackage[psamsfonts]{amssymb}
\usepackage{amsthm}
\usepackage{latexsym,xypic,euscript,array,enumerate}
\usepackage{xy}
\usepackage[dvips]{graphicx}
\usepackage[dvips]{color}
\usepackage{lscape}
\usepackage{longtable}

\numberwithin{equation}{section}

\DeclareMathOperator{\Gal}{Gal}

\newcommand{\FF}{\mathbb{F}}

\newcommand{\ZZ}{\mathbb{Z}}

\newcommand{\f}{\mathfrak{f}}

\renewcommand{\emph}[1]{\textsl{#1}}

\newcommand{\fr}{\mathfrak}

\newcommand{\tabl}[1]{\parbox{10em}{$#1$}}
\newcommand{\tabll}[1]{\parbox{27em}{$#1$}}

\theoremstyle{plain}
\newtheorem{theorem}{Theorem}[section]
\newtheorem{proposition}[theorem]{Proposition}

\newtheorem{corollary}[theorem]{Corollary}

\newtheorem{exemple}[theorem]{Example}

\newtheorem{remarks}[theorem]{Remarks}

\begin{document}
\title{$2$-extensions with many points}
\author{St{\'e}phan S{\'e}mirat}
\date{}
\begin{abstract}
We give defining equations for function fields over finite fields with many rational places. They are obtained from composita of quadratic extensions of the rational function field.
\end{abstract}
\maketitle

\section{Introduction}
\nocite{MR89b:11094}
We are interested in the following problem:
given a prime power $q=p^e$ and a positive integer $g$, what is the maximal number of rational places $N_q(g)$ of a function field of genus $g$ and finite constant field $\FF_q$? This problem arise for instance in the theory of error correcting codes.

From the Riemann Hypothesis, we have the famous Hasse-Weil bound
$$N_q(g)\leq q+1+2g\sqrt{q},$$ and J.P. Serre proved
$$N_q(g)\leq q+1+g[2\sqrt{q}],$$ 
where $[x]$ stands for the integer part of a real number $x$.

In fact, the Hasse-Weil bound is a particular case of the so-called ''\textit{m{\'e}thode des formules explicites}``, which majorates $g$ for a fixed number of rational places. However, this method doesn't give the exact value of $N_q(g)$ and to this day, there is no method which could do so. We are thus lead to trying to reach or coming as close as possible to the best known upper bounds for $N_q(g)$, by way of finding (explicitely or not) function fields with the required parameters. The tables of G. van der Geer and M. van der Vlugt \cite{Gee-Vlu5} record results on best known function fields for $1\leq g\leq 50$ and $q$ a small power of $2$ or $3$.

A classical method to obtain function fields with many rational places is given by class field theory. In this method, one has to manage with three parameters: a function field $K/\FF_q$, a set $S$ of rational places of $K$, and an effective divisor $\fr{m}$ with support disjoint from $S$. Then the ray class field $K_S^{\fr{m}}$ associated with $(K,S,\fr{m})$ has its degree over $K$ controlled by $S$ (if $S$ is sufficently large, one has $K_S^{m}=K$). The ramification of $K^{\fr{m}}_S/K$ (and thus the genus) is controlled by $\fr{m}$, and the number of rational places of $K_S^{\fr{m}}$ is greater than $[K_S^{\fr{m}}:K]\times |S|$, since all the places of $S$  split totally in $K_S^{\fr{m}}$. In particular, for $S$ large and $\deg \fr{m}$ small one can hope to obtain some ``good'' function fields.
For instance, when $\fr{m}=0$, then $K_S^{\fr{m}}$ is the Hilbert class field associated with $(K,S)$ and if $S$ equals the set of rational places of $K$, then the Hurwitz genus Formula gives
$$\frac{N(K_S^0)}{g_{K_S^0}-1}\geq \frac{N(K)}{g_K-1},$$
where we denote by $N(K)$ and $g_K$ (resp. $N(K_S^0)$ and $g_{K_S^0}$ the number of rational places and genus of $K$ (resp. $K_S^0$).\\
One may verify in the tables \cite{Gee-Vlu5} that in general, when $q$ is fixed, the quotient $N_q(g)/g$ decreases. That is the greater $N(K)$ is, the smaller $h_S(K)$ should be, when $S$ is chosen to be the whole set of rational places of $K$.\\
However, two problems arise with this method. The first one is that we don't know how to calculate the class number of a general pair $(K,S)$. The second one is that for the applications, one needs explicit equations for the function fields involved. Such equations can hardly be found using general class field theory.
\nocite{Aue2}

\section{The method}

In the following, we will apply the above class field method in it simplest form. We will construct some $(\ZZ/2\ZZ)^n$-sub-fields of some ray class field associated with quadratic extensions of the rational function field $\FF_q(x)$. The good thing with this is that we are able to give explicit equations for such fields.

Let us start with a simple example.
\begin{exemple}
Let $k=\FF_3(x)$ and let $K=k(y)$ be the quadratic extension of $k$ defined by 
\begin{equation*}
  y^2=f(x)=2(x^3+2x+2)(x^3+2x+1).
\end{equation*}
The genus of $K/\FF_3$ is $g_K=2$. Now consider the compositum $L=K_1.K_2$ where $K_1=k(y_1)$ and $K_2=k(y_2)$ are defined by the equations
$$y_1^2=f_1(x)=2(x^3+2x+2) \text{ and }y_2^2=f_2(x)=x^3+2x+1.$$
Since $f_1f_2=f$, $K$ is a subfield of $L$ and the ramification in $L/K$ occurs only at the place of degree $2$ above $(\frac{1}{x})$. The Hurwitz Formula gives then $g_L=4$, for $g_L$ the genus of $L$. Furthermore, since $f_i(0)$, $f_i(1)$, and $f_i(2)$ are all squares of $\FF_3$, we have that $N(L)=4\times 3=12$, which is the maximal number of rational places for a function fields of genus $4$ over $\FF_3$.
\end{exemple}

Now we give more general facts.\\
With the help of the genus formulae for quadratic extensions (see for instance \cite[Propositions VI.2.3 and VI.4.1]{Stibook}, or \cite{Mac}), the genus computation for a ramified $2$-elementary extension of the rational function field can be done using the following theorem, which is a particular case of a theorem of G. Frey and H. G. Rück. A proof can be found in \cite{MR97h:94022}.

\begin{theorem}
Let $L/k$ be a Galois extension such that $\Gal(L/k)\simeq (\ZZ/2\ZZ)^n$ and
let $K_I$, $ I\subset\{1,\ldots,n\}$ be the $2^n-1$ sub-fields of $L$, quadratic over $k$, and let $g_{K_I}$ be their respective genus. Then the genus $g_L$ of $L$ is given by
\begin{equation*}
g_L=\sum_ig(K_I)
\end{equation*}
\end{theorem}

\begin{remarks}
\begin{itemize}
\item In odd characteristic, we easily have $$g_L=2^{n-1}\left(\sum_{P,\ e_P>1}\deg P-4\right)+1.$$
\item An analog formula is proven in \cite{MR1732227} in the hard case of Artin-Schreier extension.
\end{itemize}
\end{remarks}

For the applications, for insance to the geometric Goppa codes, one needs basis for some Riemann-Roch vector spaces. Such basis can be obtained as soon as one possesses a defining equation for the extension $L/k$ (for instance, the software Kash \cite{Kash} knows how to do it). For our extensions $L/k$, the following proposition gives us such an equation. Therefore, we are in the best situation for the applications.

\begin{proposition}
  Let $L/k$ be the compositum of $n$ disjoint quadratic extensions $K_i/k$ given by $K_i=\FF_q(x,y_i)$ with, for each $i$, $y_i$ satisfying:
\begin{itemize}
\item in odd characteristic, $y_i^2=f_i(x)$,  where $f_i(x)\in\FF_q[x]$ is square-free,
\item in even characteristic, $y_i^2+y_i=f_i(x)$, where $f_i(x)\in\FF_q(x)$ such that for all $z\in\FF_q(x)$, $f_i(x)\neq z^2+z$.
\end{itemize}
Then $$y=\sum_iy_i\text{ (resp. $y=\prod_iy_i$)}$$ is a primitive element for $L/k$ in odd (resp. even) characteristic.
\end{proposition}
\begin{proof}
  The Galois group of $L/k$ is generated by the $\sigma_i$, $1\leq i\leq n$, such that in odd characteristic (resp. even) $\sigma_i(y_i)=-y_i$ (resp. $\sigma_i(y_i)=y_i+1$), and, for $i\neq j$, $\sigma_i(y_j)=y_j$. The element $$y=\sum_{i=1}^ny_i\text{ (resp. $y=\prod_{i=1}^n y_i$)},$$ is fixed by no  element of $\Gal(L/k)$. Thus, $L=k(y)$.\\
\end{proof}
\begin{corollary}
  Let $L=k(y_1,\ldots,y_n)=k(y)$, with $y$ defined as in the above proposition. Then a minimal polynomial for $y$ over $k$ is given by 
$$P(Y,x)=\prod_{I}\big(Y-\sigma_I(y)\big)\in k[Y],$$
where for each $I\subseteq\{1,\ldots,n\}$, 
$$\sigma_I(y)=\prod_{i=1..n,\ i\notin I}y_i\prod_{i=1..n,\ i\in I}
\begin{cases}
-y_i&\text{ in odd characteristic,}\\
y_i+1&\text{ in even characteristic.}
\end{cases}
$$
\end{corollary}
\begin{exemple}
  \begin{itemize}
  \item In odd characteristic, 
\begin{itemize}
\item the field $L=\FF_q(x,y_1,y_2)$ with $y_i^2=f_i(x)$ is given by $L=\FF_q(x,y)$ with $y$ satisfying
$$y^4-2(f_1+f_2)y^2+(f_1-f_2)^2=0,$$
\item the field $L=\FF_q(x,y_1,y_2,y_3)$ with $y_i^2=f_i(x)$ is given by $L=\FF_q(x,y)$ with $y$ satisfying\footnotemark[1]\footnotetext[1]{The subscripts in the sums aren't ordered, for instance with $1\leq i\leq 3$, the sum $\sum_{i\neq j}f_if_j$ possesses 3 terms.}
\begin{multline*} 
y^8-4\left(\sum_if_i\right)y^6+2\left(3\sum_if_i^2+2\sum_{i\neq j}f_ff_j\right)y^4-4\left(\sum_if_i^3-\sum_{i\neq j}f_i^2f_j+10f_1f_2f_3\right)y^2\\ +\sum_if_i^4-4\sum_{i\neq j}f_i^3f_j+6\sum_{i\neq j}f_i^2f_j^2+4\sum_{\substack{i\neq j,\ i\neq l\\ j\neq l}}f_i^2f_jf_l=0.
\end{multline*}
\end{itemize}
\item In even characteristic,
\begin{itemize}
\item the field $L=\FF_q(x,y_1,y_2)$ with $y_i^2+y_i^2=f_i(x)$ is given by $L=\FF_q(x,y)$ with $y$ satisfying
$$y^4+y^3+(f_1+f_2)y^2+f_1f_2y+f_1^2f_2^2=0,$$
\item the field $L=\FF_q(x,y_1,y_2,y_3)$ with $y_i^2+y_i=f_i(x)$ is given $L=\FF_q(x,y)$ with $y$ satisfying\footnotemark[1]
\begin{multline*} 
y^8+y^7+\left(\sum_if_i\right)y^6+\left(f_1f_2f_3+\sum_{i\neq j}f_if_j\right)y^5+\left(f_1f_2f_3+\sum_{i\neq j}f_i^2f_j^2\right)y^4\\
+\left(f_1^2f_2^2f_3^2+\sum_{\substack{i\neq j,\ i\neq l\\ j\neq l}}f_if_j^2f_l^2\right)y^3+\left(\sum_{\substack{i\neq j,\ i\neq l\\ j\neq l}}f_i^3f_j^2f_l^2\right)y^2+f_1^3f_2^3f_3^3y+f_1^4f_2^4f_3^4.
\end{multline*}
\end{itemize}
\end{itemize}
\end{exemple}

\section{Tables}

In the following tables, $w$ denotes a generator for $\mathbb{F}_q^*$. We give in the last column either the exact value for $N_q(g)$, or a range $[a,b]$ (written $a-b$, following the notations in \cite{Gee-Vlu5}) such that $N_q(g)\in[a,b]$ and there exists a function field defined over $\mathbb{F}_q$, with genus $g$ and with at least $a$ rational places.\\
We give the defining equations for function fields $L/\FF_q$ with genus $g_L$ obtained using the above method, and whose number of rational places $N_(L)$ is near to, or reach, or is better than $a$. In the last case $N(L)>a$, estabilhing a new minoration for the corresponding value of $N_q(g)$, we use bold face.\\
The calculation have been made with the help of the algebraic number theory software Kash, v2.2 \cite{Kash}.

\begin{landscape}
\small
\begin{longtable}{|l||c|l|c|c|c|}
\cline{1-5}
$\mathbf{q}$&$\mathbf{g_L}$ &$\mathbf{f_i}$&$\mathbf{L=\FF_q(x,y)}$  & $\mathbf{N(L)}$ & \multicolumn{1}{c}{}\\
\cline{0-5}
\endhead
\hline
$2$&$1$&\tabl{f_1=\dfrac{1}{x}\ f_2=\dfrac{1}{x+1}}&\tabll{y^4+y^3+\dfrac{1}{x(x+1)}y^2+\dfrac{1}{x(x+1)}y+\dfrac{1}{x^2(x+1)^2}}&$4$&$5$\\
\cline{2-6}
&$2$&\tabl{f_1=x\ f_2=x^3}&\tabll{y^4+y^3+(x^3+x)y^2+x^4y+x^8}&$5$&$6$\\
\cline{2-6}
&$3$&\tabl{f_1=\dfrac{1}{x}\\ f_2=\dfrac{x}{x^2+x+1}}&\tabll{y^4+y^3+\dfrac{(x+1)^2}{x(x^2+x+1)}y^2+\dfrac{1}{x(x^2+x+1)}y+\dfrac{1}{(x^2+x+1)^2}}&$6$&$7$\\
\cline{2-6}
&$4$&\tabl{f_1=x^3+x\\f^2=x+1/x}&\tabll{y^4 + y^3 + (x^4 + 1)/xy^2 + (x^4 + 1)y + x^8 + 1}&$7$&$8$\\
\cline{2-6}
&$5$&\tabl{f_1=x^3+x\\ f_2=x^5+x}&\tabll{y^4 + y^3 + (x^5 + x^3)y^2 + (x^8 + x^6 + x^4 + x^2)y + x^{16} + x^{12} + x^8 + x^4}&$9$&$9$\\
\cline{2-6}
&$7$&\tabl{f_1=x^3+x\\ f_2=\dfrac{x(x+1)}{x^3+x+1}}&\tabll{y^4+y^3+\dfrac{x^3(x+1)(x^2+x+1)}{x^3+x+1}y^2+\dfrac{x^2(x+1)^3}{x^3+x+1}y+\dfrac{x^4(x+1)^6}{(x^3+x+1)^2}}&$10$&$10$\\
\cline{2-6}
&$9$&\tabl{f_1=\dfrac{x(x+1)}{x^3+x+1}\\ f_2=\dfrac{x(x+1)}{x^3+x^2+1}}&\tabll{y^4+y^3+\dfrac{x^2(x+1)^2}{(x^3+x+1)(x^3+x^2+1)}y^2+\dfrac{x^2(x+1)^2}{(x^3+x+1)(x^3+x^2+1)}y+\dfrac{x^4(x+1)^4}{(x^3+x+1)^2(x^3+x^2+1)^2}}&$12$&$12$\\
\cline{2-6}
&$11$&\tabl{f_1=\dfrac{x(x+1)}{x^3+x+1}\\ f_2=\dfrac{x(x+1)}{x^4+x+1}}&\tabll{y^4+y^3+\dfrac{x^2(x+1)^2}{(x^3+x+1)(x^4+x+1)}y^2+\dfrac{x^2(x+1)^2}{(x^3+x+1)(x^4+x+1)}y+\dfrac{x^4(x+1)^4}{(x^3+x+1)^2(x^4+x+1)^2}}&$12$&$14$\\
\cline{2-6}
&$15$&\tabl{f_1=x+\dfrac{x}{x^2+x+1}\\ f_2=x+1+\dfrac{x+1}{x^2+x+1}\\ f_3=x^3+x}&\tabll{y^8+y^7+\dfrac{x^4(x+1)}{x^2+x+1}y^6+\dfrac{x^2(x+1)^7}{(x^2+x+1)^2}y^5+\dfrac{x^5(x+1)^5(x^4+x+1)}{(x^2+x+1)^4}y^4+\dfrac{x^6(x+1)^{12}}{(x^2+x+1)^4}y^3+\dfrac{x^{12}(x+1)^{11}}{(x^2+x+1)^5}y^2+\dfrac{x^{12}(x+1)^{15}}{(x^2+x+1)^6}y+\dfrac{x^{16}(x+1)^{20}}{(x^2+x+1)^8}}&$16$&$17$\\
\cline{2-6}
&$33$&\tabl{f_1=\dfrac{x(x+1)}{x^3+x+1}\\ f_2=\dfrac{x(x+1)}{x^3+x^2+1}\\ f_3=\dfrac{x(x+1)}{x^4+x+1}}&\tabll{y^8+y^7+\dfrac{x(x+1)(x^2+x+1)}{(x^3+x+1)(x^3+x^2+1)(x^4+x+1)}y^6+\dfrac{x^2(x+1)^2}{(x^3+x+1)(x^3+x^2+1)}y^5+\dfrac{x^3(x+1)^3(x^4+x^3+1)(x^4+x^3+x^2+x+1)}{(x^3+x+1)(x^3+x^2+1)(x^4+x+1)}y^4+(x^3+x+1)^2(x^3+x^2+1)^2(x^4+x+1)y^3+\dfrac{x^7(x+1)^7(x^2+x+1)^2}{(x^3+x+1)^3(x^3+x^2+1)^3(x^4+x+1)^3}y^2+\dfrac{x^9(x+1)^9}{(x^3+x+1)^3(x^3+x^2+1)^3(x^4+x+1)^3}y+\dfrac{x^{12}(x+1)^{12}}{(x^3+x+1)^4(x^3+x^2+1)^4(x^4+x+1)^4}}&$24$&$28-29$\\
\cline{2-6}
&$41$&\tabl{f_1=x+\dfrac{x}{x^2+x+1}\\ f_2=x+1+\dfrac{x+1}{x^2+x+1}\\ f_3=x^3+x\\ f_4=\dfrac{x(x+1)}{x^3+x+1}}&&$32$&$33-35$\\

\hline
$4$&$1$ &\tabl{f_1=x,\ f_2=\dfrac{1}{x}} &\tabll{y^4+y^3+\dfrac{(x+1)^2}{x}y^2+y+1}& $8$ & $9$ \\
\cline{2-6}
&$2$ &\tabl{ f_1=x,\ f_2=x^3} &\tabll{y^4+y^3+(x^3+x)y^2+x^4y+x^8}&$ 9$ &$ 10$ \\
\cline{2-6}
&$3$ &\tabl{ f_1=\dfrac{1}{x}+\dfrac{w}{x+w^2}\\ f_2=\dfrac{1}{x}+\dfrac{w^2}{x+w}} &\tabll{y^4+y^3+\dfrac{x+1}{(x+w)(x+w^2)}y^2+\dfrac{(x+1)^2}{x^2(x+w)(x+w^2)}y+\dfrac{(x+1)^4}{x^4(x+w)^2(x+w^2)^2}}&$ 14$ &$ 14$ \\
\cline{2-6}
&$4$&\tabl{f_1=x^3+1\ f_2=x+1/x}&\tabll{y^4+y^3+\dfrac{(x+1)(x^3+x^2+1)}{x}y^2+\dfrac{(x+1)^3(x+w)(x+w^2)}{x}y+(x+1)^6(x+w)^2(x+w^2)^2}&$15$&$15$\\
\cline{2-6}
&$5$ &\tabl{ f_1=x^3+1,\ f_2=x^5+x}&\tabll{y^4+y^3+(x+1)(x^2+wx+w)(x^2+w^2x+w^2)y^2+x(x+1)^5(x+w)(x+w^2)y+x^2(x+1)^{10}(x+w)^2(x+w^2)^2}&$ 17$ &$ 17-18$ \\
\cline{2-6}
&$7$ &\tabl{f_1=x^3+1\\ f_2=\dfrac{x(x+w^2)}{x^3+w^2}}&\tabll{y^4 + y^3 + (x^6 + wx^3 + x^2 + w^2x + w^2)/(x^3 + w^2)y^2 + (x^5 + w^2x^4 + x^2 + w^2x)/(x^3 + w^2)y + (x^{10} + wx^8 + x^4 + wx^2)/(x^6 + w)}&$ 18$ &$ 21-22$ \\
\cline{2-6}
&$11$ &\tabl{ f_1=x^3+1\\ f_2=x+\dfrac{1}{x}\\ f_3=x+1+\dfrac{1}{x+1}} &\tabll{y^8+y^7+\dfrac{(x+w)(x+w^2)(x^3+x+1)}{x(x+1)}y^6+\dfrac{(x^3+wx^2+1)(x^3+w^2x^2+1)}{x}y^5+\dfrac{(x^3+x^2+1)(x^5+x^2+1)}{x^2}y^4+(x+1)^2(x+w)^2(x+w^2)(x^3+wx^2+1)(x^3+w^2x^2+1)y^3+x(x+1)^3(x+w)^3(x+w^2)^3(x^3+x+1)y^2+xx^3(x+1)^6(x+w)^3(x+w^2)^3y+x^4(x+1)^8(x+w)^4(x+w^2)^4}&$ 26$ &$ 26-30$ \\
\cline{2-6}
\hline
\hline
$8$&$1$ &\tabl{f_1=\dfrac{1}{x},\ f_2=\dfrac{w}{x+w^4}} &\tabll{y^4 + y^3 + (w^3x + w^4)/(x^2 + w^4x)y^2 + w/(x^2 + w^4x)y + w^2/(x^4 + wx^2)}&$ 14 $&$ 14$ \\
\cline{2-6}
&$2$ &\tabl{ f_1=\dfrac{1}{x}+\dfrac{w^6}{x+1}\\  f_2=\dfrac{w^3}{x}} &\tabll{y^4+y^3+\dfrac{w^5(x+w^3)}{x(x+1)}y^2+\dfrac{w^5(x+w^5)}{x^2(x+1)}y+\dfrac{w^3(x+w^5)^2}{x^4(x+1)^2}}&$ 17$ &$ 18$ \\
\cline{2-6}
&$3$ &\tabl{ f_1=\dfrac{1}{x}+\dfrac{w}{x+w^4},\ f_2=\dfrac{1}{x}+\dfrac{w^5}{x+w^6}}&\tabll{y^4+y^3+\dfrac{w^6(x+1)}{(x+w^4)(x+w^6)}y^2+\dfrac{(x+w)(x+w^2)}{x^2(x+w^4)(x+w^6)}y+\dfrac{(x+w)^2(x+w^2)^2}{x^4(x+w^4)(x+w^6)}}& $24$ &$ 24$ \\
\cline{2-6}
&$4$ &\tabl{ f_1=x^3+x+1\\ f_2=x+\dfrac{1}{x+1}}&\tabll{y^4+y^3+\dfrac{x(x+w^3)(x+w^5)(x+w^6)}{x+1}y^2+\dfrac{(x+w)(x+w^2)(x+w^4)(x^2+x+1)}{x+1}y+\dfrac{(x+w)^2(x+w^2)^2(x+w^4)^2(x^2+x+1)^2}{(x+1)^2}}& $25$ &$ 25-27$ \\
\cline{2-6}
&$5$ &\tabl{ f_1=\dfrac{1}{x^2+x+1}\\ f_2=\dfrac{w(x+w^6)}{x^2+w^6x+w}} &\tabll{y^4 + y^3 + (wx^3 + wx^2 + w^4x + w^3)/(x^4 + w^2x^3 + w^4x^2 + w^5x + w)y^2 + (wx + 1)/(x^4 + w^2x^3 + w^4x^2 + w^5x + w)y + (w^2x^2 + 1)/(x^8 + w^4x^6 + wx^4 + w^3x^2 + w^2)}&$ 28$ &$ 29-32$ \\
\cline{2-6}
&$7$ &\tabl{f_1=\dfrac{1}{x}+\dfrac{w(x+w^3)}{x^2+w^5x+w}\\ f_2=\dfrac{1}{x}+\dfrac{w^2(x+w^{6})}{x^2+w^3x+w^2}} &\tabll{y^4 + y^3 + (w^4x^3 + w^3x^2 + w^2x + 1)/(x^4 + w^2x^3 + w^2x^2 + w^5x + w^3)y^2 + (w^2x^4 + w^4x^3 + w^5x^2 + w^4x + w^3)/(x^6 + w^2x^5 + w^2x^4 + w^5x^3 + w^3x^2)y + (w^4x^8 + wx^6 + w^3x^4 + wx^2 + w^6)/(x^{12} + w^4x^{10} + w^4x^8 + w^3x^6 + w^6x^4)}&$ \textbf{34}$ &$ 33-39$ \\
\cline{2-6}
&$8$&\tabl{ f_1=1/x+1/(x + w^2)\\ f_2=w/(x+w^2)+ w^6/(x + w^6)\\f_3=w^5/(x+w^2)}&\tabll{y^8 + y^7 + (x + w)/(x^3 + x^2 + wx)y^6 + (w^3x^2 + w^6x + w^4)/(x^4 + x^3 + w^2x^2 + w^4x + w^5)y^5 + (wx^3 + w^5x^2 + w^3x + 1)/(x^7 + w^2x^6 + x^5 + w^2x^4 + w^2x^3 + w^4x^2)y^4 + (wx^3 + w^3x^2 + 1)/(x^9 + x^7 + w^4x^5 + wx^3 + w^3x)y^3 + (w^3x^3 + w^4x^2 + w^6x + 1)/(x^{13} + x^{12} + w^3x^{11} + x^{10} + w^2x^9 + w^4x^8 + w^6x^7 + wx^6 + w^5x^5 + w^3x^4 + w^4x^3)y^2 + (wx^3 + w^6x^2 + w^4x + w^2)/(x^{15} + x^{14} + w^6x^{13} + w^5x^{12} + w^6x^{11} + w^2x^7 + w^2x^6 + wx^5 + x^4 + wx^3)*y + (w^6x^4 + w^5)/(x^{20} + x^{16} + wx^{12} + w^2x^8 + w^6x^4)}&$33$&$34-43$\\
\cline{2-6}
&$11$ &\tabl{f_1=x^3+x+1\\ f_2=x+\dfrac{1}{x+1}\\ f_3=x+1+\dfrac{1}{x}}&\tabll{y^8 + y^7 + (x^5 + x^4 + x^3 + x^2 + 1)/(x^2 + x)y^6 + (x^7 + x^5 + x^4 + x + 1)/(x^2 + x)y^5 + (x^8 + x^6 + x^5 + x^3 + x^2 + x + 1)/(x^2 + x)y^4 + (x^{14} + x^{12} + x^8 + x^7 + x^5 + x^3 + 1)/(x^4 + x^2)y^3 + (x^{19} + x^{18} + x^{17} + x^{16} + x^{14} + x^{13} + x^{12} + x^{11} + x^{10} + x^9 + x^4 + x^3 + 1)/(x^6 + x^5 + x^4 + x^3)y^2 + (x^{21} + x^{18} + x^{16} + x^{14} + x^{12} + x^{11} + x^{10} + x^7 + x^5 + x^4 + x^3 + x + 1)/(x^6 + x^5 + x^4 + x^3)y + (x^{28} + x^{16} + x^8 + x^4 + 1)/(x^8 + x^4)}&$48$&$48-54$\\
\hline
\hline
$16$&$1$ &\tabl{f_1=\dfrac{1}{x},\ f_2=\dfrac{1}{x+w^5}} &\tabll{y^4 + y^3 + w^5/(x^2 + w^5x)y^2 + 1/(x^2 + w^5x)y + 1/(x^4 + w^{10}x^2)}&$ 24$ &$ 25$ \\
\cline{2-6}
&$2$ &\tabl{f_1=x^3+x+1\\ f_2=w^5x} &\tabll{y^4 + y^3 + (x^3 + w^{10}x + 1)y^2 + (w^5x^4 + w^5x^2 + w^5x)y + w^{10}x^8 + w^{10}x^4 + w^{10}x^2}&$ 33$ &$ 33$ \\
\cline{2-6}
&$3$ &\tabl{f_1=\dfrac{1}{x}+\dfrac{1}{x+w^{10}}\\ f_2=\dfrac{1}{x}+\dfrac{w^5}{x+w^5}}&\tabll{y^4 + y^3 + (w^{10}x + w^{10})/(x^2 + x + 1)y^2 + w^5/(x^3 + w^5x^2)y + w^{10}/(x^6 + w^{10}x^4)}&$ 38 $& $38$ \\
\cline{2-6}
&$4$&\tabl{f_1:=\frac{w^7}{x}+\frac{w}{x+w^{14}}\\ f_2=\frac{w^2}{x}+\frac{1}{x+w^6}}&\tabll{y^4 + y^3 + (w^6x^2 + w^2x + w^2)/(x^3 + w^8x^2 + w^5x)y^2 + (w^7x^2 + wx + w^{14})/(x^4 + w^8x^3 + w^5x^2)y + (w^{14}x^4 + w^2x^2 + w^{13})/(x^8 + wx^6 + w^{10}x^4)}&$45$&$45-46$\\
\cline{2-6}
&$5$ &\tabl{f_1=x^3+x\\ f_2=x^5}&\tabll{y^4 + y^3 + (x^5 + x^3 + x)y^2 + (x^8 + x^6)y + x^{16} + x^{12}}&$ 49$ &$ 49-54$ \\

\hline
\hline
$32$&$1$ &\tabl{f_1=\dfrac{1}{x},\ f_2=\dfrac{1}{x+1}} &\tabll{y^4 + y^3 + y^2/(x^2 + x) + y/(x^2 + x) + 1/(x^4 + x^2)}&$ 44 $&$ 44$ \\
\cline{2-6}
&$3 $&\tabl{f_1=\dfrac{1}{x}+\dfrac{1}{x+1}\\ f_2=\dfrac{1}{x}+\dfrac{w}{x+w^9}}&\tabll{y^4+y^3+\dfrac{w^{18}(x+w^3)}{(x+1)(x+w^9)}y^2+\dfrac{w^{18}(x+w^{22})}{x^2(x+1)(x+w^9)}+\dfrac{w^5(x+w^{22})}{x^4(x+1)^2(x+w^9)^2}}&$ \textbf{64} $&$ 63-64 $\\
\cline{2-6}
&$5 $&\tabl{f_1=\dfrac{1}{x^2+w^3x+w}\\ f_2=\dfrac{x+w^{24}}{x^2+w^{27}x+w^{10}}}&\tabll{y^4 + y^3 + (x^3 + w^{26}x^2 + wx + w^3)/(x^4 + w^{18}x^3 + x^2 + w^6x + w^{11})y^2 + (x + w^{24})/(x^4 + w^{18}x^3 + x^2 + w^6x + w^{11})y + (x^2 + w^{17})/(x^8 + w^5x^6 + x^4 + w^{12}x^2 + w^{22})}&$76$&$83-86$\\
\cline{2-6}
&$7 $&\tabl{f_1=\dfrac{1}{x}+\dfrac{1}{x+1}\\ f_2=\dfrac{1}{x}+\dfrac{w^{11}}{x+w^{21}}\\ f_3=\dfrac{w^{21}}{x}}&\tabll{y^8 + y^7 + (w^{24}x^2 + w^{11})/(x^3 + w^{25}x^2 + w^{21}x)y^6 + (w^{9}x^3 + w^{18}x^2 + wx + w^{11})/(x^5 + w^{25}x^4 + w^{21}x^3)*y^5 + (w^{25}x^4 + w^{23}x^3 + w^{30}x^2 + wx + w^{11})/(x^8 + w^{19}x^6 + w^{11}x^4)y^4 + (w^{18}x^4 + w^{11}x^3 + w^{21}x^2 + wx + w^{22})/(x^{10} + w^{19}x^8 + w^{11}x^6)y^3 + (w^{11}x^4 + w^{28}x^2 + w^{2})/(x^{13} + w^{25}x^{12} + w^{24}x^{11} + w^{13}x^{10} + w^{14}x^9 + w^{5}x^8 + wx^7)y^2 + (w^{27}x^3 + w^{29}x^2 + x + w^{2})/(x^{15} + w^{25}x^{14} + w^{24}x^{13} + w^{13}x^{12} + w^{14}x^{11} + w^{5}x^{10} + wx^9)y + (w^{5}x^4 + w^{13})/(x^{20} + w^{7}x^{16} + w^{22}x^{12})}&\textbf{98}&$ 90-108 $\\
\cline{2-6}
&$9 $&\tabl{f_1=\dfrac{1}{x}+\dfrac{1}{x+1}\\ f_2=\dfrac{1}{x}+\dfrac{w}{x+w^9}\\ f_3=\dfrac{w^{24}}{x+w^5}}&\tabll{y^8 + y^7 + (w^{14}x^2 + w^8x + w^{17})/(x^3 + w^{24}x^2 + wx + w^{14})y^6+ (w^{11}x^3 + w^{24}x^2 + w^{30}x + w^{25})/(x^5 + w^{24}x^4 + wx^3 + w^{14}x^2)y^5 + (w^{30}x^6 + w^7x^5 + w^6x^4 + w^{10}x^3 + w^{11}x^2 + w^{28})/(x^{10} + w^{17}x^8 + w^2x^6 + w^{28}x^4)y^4 + (w^{22}x^4 + w^{20}x^3 + w^{19}x^2 + w^{11}x + w^{27})/(x^{10} + w^{17}x^8 + w^2x^6 + w^{28}x^4)y^3 + (w^5x^4 + w^{30}x^3 + w^{12}x^2 + w^{12}x + w^{21})/(x^{13} + w^{24}x^{12} + w^{10}x^{11}+ w^{20}x^{10} + w^{11}x^9 + w^{28}x^8 + w^{24}x^7 +w^{18}x^6 + w^{29}x^5 + w^{11}x^4)y^2 + (w^2x^3 + w^{24}x^2 + w{^15}x + w^6)/(x^{15} + w^{24}x^{14} + w^{10}x^{13} + w^{20}x^{12} + w^{11}x^{11} + w^{28}x^{10} + w^{24}x^9 + w^{18}x^8 + w^{29}x^7 + w^{11}x^6)y + (w^{13}x^4 + w^8)/(x^{20} + w^3x^{16} + w^4x^{12} + w^{25}x^8)
} &$ 104 $&$ 108-130 $\\
\cline{2-6}
&$11$&\tabl{f_1=\dfrac{1}{x}+\dfrac{1}{x+1}\\ f_2=\dfrac{1}{x}+\dfrac{w^{4}}{x+w^{2}}\\ f_3=\dfrac{w^7}{x}+\dfrac{w^{13}}{x+w^{22}}}&\tabll{y^8 + y^7 + (w^{25}x^3 + w^{30}x^2 + w^{2}x + 1)/(x^4 + w^{4}x^3 + w^{23}x^2 + w^{24}x)y^6 + (w^{13}x^4 + w^{8}x^3 + w^{15}x^2 + w^{28}x + 1)/(x^6 + w^{4}x^5 + w^{23}x^4 + w^{24}x^3)y^5 + (w^{27}x^6 + w^{16}x^5 + w^{27}x^4 + w^{8}x^3 + x^2 + w^{24}x + w^{17})/(x^{10} + w^{8}x^8 + w^{15}x^6 + w^{17}x^4)y^4 + (w^{26}x^6 + w^{18}x^5 + w^{8}x^4 + w^{15}x^3 + wx^2 + w^{24}x + 1)/(x^{12} + w^{8}x^{10} + w^{15}x^8 + w^{17}x^6)y^3 + (w^{20}x^7 + w^{25}x^6 + w^{26}x^5 + w^{30}x^4 + w^{7}x^3 + w^{21}x^2 + w^{2}x + 1)/(x^{16} + w^{4}x^{15} + wx^{14} + w^{4}x^{13} + w^{24}x^{12} + w^{2}x^{11} + w^{11}x^{10} + w^{22}x^9 + w^{9}x^8 + w^{10}x^7)y^2 + (w^{8}x^6 + w^{29}x^5 + w^{10}x^4 + w^{9}x^3 + w^{28}x^2 + w^{3}x + 1)/(x^{18} + w^{4}x^{17} + wx^{16} + w^{4}x^{15} + w^{24}x^{14} + w^{2}x^{13} + w^{11}x^{12} + w^{22}x^{11} + w^{9}x^{10} + w^{10}x^9)y + (w^{21}x^8 + w^{12}x^4 + 1)/(x^{24} + w^{16}x^{20} + w^{30}x^{16} + w^{3}x^{12})}&\textbf{120}&$113-152$\\
\hline
\hline
$64$&$1 $&\tabl{f_1=\dfrac{1}{x},\ f_2=\dfrac{1}{x+w^{27}}}&\tabll{y^4 + y^3 + w^{27}/(x^2 + w^{27}x)y^2 + 1/(x^2 + w^{27}x)y + 1/(x^4 + w^{54}x^2)}&$ 80 $&$ 81 $\\
\cline{2-6}
&$3 $&\tabl{f_1=\dfrac{1}{x^2+w^{25}x+1}\\ f_2=\dfrac{w^{14}}{x+w^{28}}} &\tabll{y^4 + y^3 + (w^{14}x^2 + w^{17}x + w^{30})/(x^3 + w^{38}x^2 + w^{50}x + w^{28})y^2 + w^{14}/(x^3 + w^{38}x^2 + w^{50}x + w^{28})y + w^{28}/(x^6 + w^{13}x^4 + w^{37}x^2 + w^{56})}&$ 104 $&$ 113 $\\
\cline{2-6}
&$5 $&\tabl{f_1=\dfrac{1}{x^2+w^{11}x+1}\\ f_2=\dfrac{1}{x^2+w^{25}x+1} }&\tabll{y^4 + y^3 + w^{27}x/(x^4 + w^{27}x^3 + w^{36}x^2 + w^{27}x + 1)y^2 + 1/(x^4 + w^{27}x^3 + w^{36}x^2 + w^{27}x + 1)y + 1/(x^8 + w^{54}x^6 + w^9x^4 + w^{54}x^2 + 1)}&$ 128 $&$ 130-145 $\\
\hline
\hline
$128$&$1 $&\tabl{f_1=\dfrac{1}{x},\ f_2=\dfrac{w^{11}}{x+w^{111}} }&\tabll{y^4 + y^3 + (w^{87}x + w^{111})/(x^2 + w^{111}x)y^2 + w^{11}/(x^2 + w^{111}x)y + w^{22}/(x^4 + w^{95}x^2)}&$ 150 $&$ 150 $\\
\cline{2-6}
&$3 $&\tabl{f_1=\dfrac{1}{x}+\dfrac{w}{x+w^{19}}\\ f_2=\dfrac{1}{x}+\dfrac{w^{74}}{x+w^{79}}}&\tabll{y^4 + y^3 + (w^{79}x + w^8)/(x^2 + w^{16}x + w^{98})y^2 + (w^{99}x^2 + w^{64}x + w^{98})/(x^4 + w^{16}x^3 + w^{98}x^2)y + (w^{71}x^4 + wx^2 + w^{69})/(x^8 + w^{32}x^6 + w^{69}x^4)}&$\textbf{192}$&$191-192$\\
\cline{2-6}
&$5 $&\tabl{f_1=\dfrac{1}{x}\\ f_2=\dfrac{w^{11}}{x+w^{111}}\\ f_3=\dfrac{w^{42}}{x+w^{63}} }&\tabll{y^8 + y^7 + (w^{110}x^2 + w^{64}x + w^{47})/(x^3 + w^{55}x^2 + w^{47}x)y^6 + (w^{92}x + 1)/(x^3 + w^{55}x^2 + w^{47}x)y^5 + (w^{53}x^3 + w^{97}x^2 + w^{100}x + w^{36})/(x^6 + w^{110}x^4 + w^{94}x^2)y^4 + (w^{18}x + w^{53})/(x^6 + w^{110}x^4 + w^{94}x^2)y^3 + (w^{89}x^2 + w^{43}x + w^{26})/(x^9 + w^{55}x^8 + w^{50}x^7 + w^{38}x^6 + w^{97}x^5 + w^{22}x^4 + w^{14}x^3)y^2 + w^{32}/(x^9 + w^{55}x^8 + w^{50}x^7 + w^{38}x^6+ w^{97}x^5 + w^{22}x^4 + w^{14}x^3)y + w^{85}/(x^{12} + w^{93}x^8 + w^{61}x^4)}&$ 216 $&$ 227-239
$\\
\cline{2-6}
&$9 $&\tabl{f_1=\dfrac{1}{x}+\dfrac{w}{x+w^{19}}\\ f_2=\dfrac{1}{x}+\dfrac{w^{74}}{x+w^{79}}\\ f_3=\dfrac{1}{x}+\dfrac{w^{80}}{x+w^{61}}}&\tabll{y^8 + y^7 + (w^{34}x^3 + w^{24}x^2 + w^{39}x + w^{32})/(x^4 + w^{84}x^3 + w^{7}x^2 + w^{32}x)y^6 + (w^{122}x^4 + w^{33}x^3 + w^{112}x^2 + w^{33}x + w^{32})/(x^6 + w^{84}x^5 + w^{7}x^4 + w^{32}x^3)y^5 + (w^{74}x^7 + w^{109}x^6 + w^{35}x^5 + w^{51}x^4 + w^{15}x^3 + w^{63}x^2 + w^{64}x + w^{64})/(x^{10} + w^{41}x^8 +w^{14}x^6 + w^{64}x^4)y^4 + (w^{69}x^7 + w^{2}x^6 + w^{30}x^5 + w^{6}x^4 + w^{47}x^3 + w^{116}x^2 + w^{64}x + w^{64})/(x^{12} + w^{41}x^{10} + w^{14}x^8 + w^{64}x^6)y^3 + (w^{55}x^9 + w^{45}x^8 + w^{118}x^7 + w^{35}x^6 + w^{66}x^5 + w^{95}x^4 + x^3 + w^{8}x^2 + w^{103}x + w^{96})/(x^{16} + w^{84}x^{15} + w^{93}x^{14} +w^{84}x^{13} + w^{100}x^{12} + w^{51}x^{11} + w^{38}x^{10} + w^{126}x^9 + w^{71}x^8 + w^{96}x^7)y^2 + (w^{95}x^9 + w^{53}x^8 + w^{6}x^7 + w^{70}x^6 + w^{98}x^5 + w^{98}x^4 + w^{47}x^3 + w^{98}x^2 + w^{103}x + w^{96})/(x^{18} + w^{84}x^{17} + w^{93}x^{16} + w^{84}x^{15} + w^{100}x^{14} + w^{51}x^{13} + w^{38}x^{12} + w^{126}x^{11} + w^{71}x^{10} + w^{96}x^9)y + (w^{42}x^{12} + wx^8 + w^{29}x^4 + w)/(x^{24} + w^{82}x^{20} + w^{28}x^{16} + wx^{12})}&$ \textbf{288} $&$ 258-327 $\\
\hline
\hline
$3$&$1 $&\tabl{f_1=x\\ f_2=x^2+x+2 }&\tabll{y^4 + (x^2 + 2x + 2)y^2 + x^4 + x^2 + 1}&$ 6 $&$ 7 $\\
\cline{2-6}
&$2 $&\tabl{f_1=x\\ f_2=x^3+2x+1}&\tabll{y^4 + (x^3 + 1)y^2 + x^6 + 2x^4 + 2x^3 + x^2 + 2x + 1}&$ 8 $&$8$\\
\cline{2-6}
&$3 $&\tabl{f_1=(x+1)(x^3+x^2+x+2)\\ f_2=x^2+x+2 }&\tabll{y^4 + (x^4 + 2x^3 + x + 1)y^2 + x^8 + x^7 + 2x^5 + x^3 + x^2}&$ 8 $&$ 10 $\\
\cline{2-6}
&$4 $&\tabl{f_1=x^3+2x+1\\ f_2=2(x^3+2x+2)}&\tabll{y^4 + 2y^2 + x^6 + x^4 + x^2}&$ 12 $&$ 12 $\\
\cline{2-6}
&$5 $&\tabl{f_1=(x+1)(x^3+x^2+2x+1)\\ f_2=x^4+x^3+x^2+1 }&\tabll{y^4 + (2x^4 + x^2 + 2)y^2 + x^6 + x^5 + x^4}&$ 12 $&$ 12-13 $\\
\cline{2-6}
&$6 $&\tabl{f_1=x^3+2x+1\\ f_2=x^5+2x+1}&\tabll{y^4 + (x^5 + x^3 + x + 2)y^2 + x^{10} + x^8 + x^6}&$ 14 $&$ 14-15 $\\
\cline{2-6}
&$7 $&\tabl{f_1=(x^2+1)(x^2+x+2)(x^2+2x+2)\\ f_2=(x^2+1)(x^4+x^3+2x+1)}&\tabll{y^4 + (2x^6 + x^5 + 2x^4 + 2x^2 + 2x + 2)y^2 + x^{10} + x^6 + x^2}&$ 16 $&$ 16 $\\
\cline{2-6}
&$8 $&\tabl{f_1=x^5+2x+1\\ f_2=x^5+x^3+x+1}&\tabll{y^4 + (2x^5 + x^3 + 2)y^2 + x^6 + x^4 + x^2}&$ 14 $&$ 15-18 $\\
\cline{2-6}
&$9 $&\tabl{f_1=x^2+1\\ f_2=(x^2+x+2)(x^2+2x+2)\\ f_3=(x+1)(x+2)(x^2+x+2)}&\tabll{y^8 + (x^4 + 2x^3 + x^2 + x)y^6 + (x^8 + x^7 + 2x^4 + x^3 + x^2 + x)*y^4 + (2x^{10} + 2x^9 + x^8 + x^7 + 2x^6 + 2x^4 + 2x^3 + x + 2)y^2 + x^6 + x^4 + x^2}&$ 16 $&$ 19 $\\
\cline{2-6}
&$13 $&\tabl{f_1=(x+1)(x^3+x^2+2x+1)\\ f_2=(x+1)(x^2+1) \\ f_3=x^3+2x+1}&\tabll{y^8 + (2x^4 + 2x^3 + 2x^2)y^6 + (2x^7 + 2x^4 + x^3 + x^2)y^4 + (2x^{12}+ 2x^{11} + x^{10} + x^9 + 2x^7 + x^5 + 2)y^2 + x^{16} + x^{13} + x^{12} + x^{11} + 2x^{10} + 2x^9 + x^7 + 2x^5 + x^4}&$ 24 $&$ 24-25 $\\
\cline{2-6}
&$17 $&\tabl{f_1=(x+2)(x^3+x^2+x+2)\\ f_2=x^4+x^3+x^2+x+1\\ f_3=x^4+2x^3+x^2+1 }&\tabll{y^8 + (x^2 + x)y^6 + (x^5 + 2x^2 + x)y^4 + (2x^{12} + 2x^{10} + x^9 + 2x^6 + 2x^5 + x^3 + 2x^2 + x + 2)y^2 + x^{10} + x^9 + x^6 + 2x^3 + x^2}&$ 24 $&$ 24-30 $\\
\cline{2-6}
&$19 $&\tabl{f_1=(x^2+1)(x^2+x+2)(x^2+2x+2)\\ f_2=(x^2+1)(x^3+x^2+2x+1)\\ f_3=x^3+2x+1}&\tabll{y^8 + (2x^6 + 2x^5 + x^4 + 2x^3 + 2x)y^6 + (x^{11} + x^{10} + 2x^9 + x^8 + x^7 + x^6 + 2x^5 + 2x^4 + x^3 + x^2 + 2x)y^4 + (2x^{18} + x^{17} + 2x^{16} + x^{15} + 2x^{14} + 2x^{12} + x^{11} + 2x^7 + x^6 + 2x^5 + 2x^4 + x^2 + 2x + 2)y^2 + x^{24} + 2x^{23} + x^{21} + x^{20} + x^{19} + x^{18} + x^{17} + x^{16} + 2x^{15} + 2x^{14} + x^{13} + 2x^{10} + x^7 + 2x^6 + x^4 + x^2}&$ 28 $&$ 28-32 $\\
\cline{2-6}
&$21 $&\tabl{f_1=(x^2+1)(x^2+x+2)(x^2+2x+2)\\ f_2=(x^2+1)(x^4+x^3+2x+1)\\ f_3=(x^2+x+2)(x^4+x+2)}&\tabll{y^8 + (x^5 + 2x^4 + 2x^3 + x^2)y^6 + (x^{11} + x^9 + x^8 + 2x^7 + 2x^5 + 2x^4)y^4 + (2x^{18} + x^{17} + 2x^{15} + x^{12} + x^{11} + 2x^9 + x^6 + 2x^5 + x^3 + 2x^2 + 2)y^2 + x^{22} + x^{21} + x^{20} + x^{19} + x^{18} + x^{16} + x^{15} + x^{12} + x^{11} + 2x^9 + x^8}&$ 32 $&$ 32-35 $\\
\hline
\hline
$9$&$1 $&\tabl{f_1=x(x+w^6)\\ f_2=(x+w^2)}&\tabll{y^4 + (x^2 + w^5x + w^2)y^2 + x^4 + w^7x^3 + w^5x + 2}&$ 16   $&$ 16 $\\
\cline{2-6}
&$2$&\tabl{f_1=x\\ f_2=(x+1)(x+w^5)(x+w^7)}&\tabll{y^4 + (x^3 + 2*x + 2)*y^2 + x^6 + x^3 + 1}&$18$&$20$\\
\cline{2-6}
&$3 $&\tabl{f_1=x(x+w)(x+1)(x+w^3)\\ f_2=x(x+w)(x+w^6)(x+w^7)}&\tabll{y^4 + (2x^4 + wx^3 + w^3x)y^2 + w^6x^6 + x^4 + w^2x^2}&$ 28 $&$ 28 $\\
\cline{2-6}
&$4 $&\tabl{f_1=x^3 + w^2x^2 + w^2\\ f_2=x^3 + w^6x^2 + w^6}&\tabll{y^4 + 2x^3y^2 + 2x^4 + x^2 + 2}&$ 30 $&$ 30 $\\
\cline{2-6}
&$5 $&\tabl{f_1=x(x+w)\\ f_2=(x+1)(x+w^3)\\ f_3=(x+w^6)(x+w^7)}&\tabll{y^8 + w^6y^6 + y^4 + (2x^6 + w^2x^4 + w^3x^3 + x^2 + w^5x + w^6)y^2 + 1}&$ 32 $&$ 32-36 $\\
\cline{2-6}
&$6 $&\tabl{f_1=x(x^2 + w)(x^3 + wx^2 + 2x + w^7)\\  f_2=x(x^2 + w)(x^3 + w^3x^2 + w^6x + w^3)}&\tabll{y^4 + (2x^6 + x^5 + w^2x^4 + wx^3 + 2x^2)y^2 + 2x^{10} + w^3x^9 + w^2x^8 + 2x^7 + x^6 + wx^5 + 2x^4 + w^2x^3 + x^2}&$ 34 $&$ 35-40 $\\
\cline{2-6}
&$7 $&\tabl{f_1=(x^2 + w^2x + w^6)(x^4 + wx^2 + wx + w^2)\\ f_2=(x^2 + w^2x + w^6)(x^4 + w^7x^2 + w^5x + w^6)}&\tabll{y^4 + (2x^6 + w^6x^5 + x^3 + 2x^2)y^2 + x^8 + 2x^7 + w^2x^6 + w^3x^5 +w^5x^4 + wx^3 + w^6x^2 + w^6x + 1}&$ 36 $&$ 40-43 $\\
\cline{2-6}
&$8 $&\tabl{f_1=(x^3+w^2x^2+2x+w^2)(x^5+2x+w^2)\\ f_2=(x^3+w^6x^2+x+w^2)(x^5+2x+w^2)}&\tabll{y^4 + (2x^8 + w^6x^5 + x^4 + w^6x^3 + w^2x + 1)y^2 + 2x^{14} + w^2x^{13} + x^{12} + 2x^{10} + w^6x^9 + 2x^8 + w^6x^7 + 2x^6 + x^4 + 2x^2}&$\textbf{40} $&$ 38-47 $\\
\cline{2-6}
&$9 $&\tabl{f_1=w*x^4 + w^7,\\ f_2=x^4 + x^3 + w^5*x + w^6 \\ f_3=x^4 + w^3*x^3 + w*x^2 + 2*x}&\tabll{y^8 + (w^3x^4 + w^2x^3 + w^5x^2 + w^2x + 2)y^6 + (w^3x^8 + x^7 + w^7x^6 + w^2x^5 + w^6x^4 + w^6x^3 + w^5x^2 + x + w^5)y^4 + (w^2x^{12} + wx^{11} + w^6x^9 + w^2x^8 + w^3x^7 + w^2x^6 + wx^5 + w^2x^4 + w^7x^3 + w^2x^2 + w^7x + w^6)y^2 + x^{16} + w^3x^{15} + w^5x^{13} + w^6x^{12} + w^6x^{11} + wx^{10}+ w^2x^9 + w^5x^8 + wx^7 + 2x^6 + wx^5 + w^3x^4 + w^6x^3 + w^6x^2 + w^7x + 2}&$ 40 $&$ 48-51 $\\
\cline{2-6}
&$17 $&\tabl{f_1=(x^3+wx^2+wx+w^6)(x^3+w^3x^2+w^3x+w^2)\\ f_2=(x^3+wx^2+wx+w^6)(x^3+x+w^2)\\ f_3=(x^3+wx^2+wx+w^6)(x^3+x+w^6)}&\tabll{y^8 + (w^7x^5 + w^2x^4 + x^3 + wx^2 + w^3x + 2)y^6 + (w^7x^{11} + w^7x^{10} + 2x^9 + wx^8 + x^5 + 2x^4 + w^7x^3 + 2x^2 + 2x + 2)y^4 + (2x^{18} + w^7x^{17} + w^3x^{16} + wx^{15} + w^3x^{14} + w^6x^{13} + w^7x^{12} + w^2x^{11} + w^6x^{10} + 2x^9 + w^5x^7 + w^7x^6 + w^7x^5 + w^2x^2 + 2)y^2 + w^6x^{22} + x^{21} + w^6x^{20} + 2x^{19} + w^6x^{18} + w^5x^{17} + w^7x^{16} + w^6x^{15} + w^3x^{14} + 2x^{13} + wx^{12} + w^6x^{11} + wx^9 + w^6x^8 + w^7x^6 + w^2x^4 + w^5x^3 + 2x^2 + w^6x + 1}&$ 64 $&$ 64-82 $\\
\hline
\hline
$27$&$1 $&\tabl{f_1=x(x+w^{14})\\ f_2=x^2+w^2x+w^{10} }&\tabll{y^4 + (2x^2 + w^4x + w^{10})y^2 + w^{20}x^2 + w^7x + w^{20}}&$ 38 $&$ 38 $\\
\cline{2-6}
&$3 $&\tabl{f_1=(x+w^5)(x^2+w^{12}x+w^3)\\ f_2=x^2+w^{18}x+w^{25}}&\tabll{y^4 + (x^3 + w^3x^2 + wx + w^2)y^2 + x^6 + w^{14}x^5 + 2x^4 + w*x^3 + w^9x^2 + w^8x + 1}&$ 52 $&$ 56 $\\
\cline{2-6}
&$4 $&\tabl{ f_1=(x+w)(x^2+w^2)\\ f_2=x^3+w^7x+w^{25}}&\tabll{y^4 + (2x^3 + wx^2 + w^{19}x + w^{17})y^2 + w^2x^4 + w^{23}x^3 + 2x^2 + w^2x + w^{12}}&$ 64 $&$ 64-66 $\\
\cline{2-6}
&$5 $&\tabl{f_1=x(x + w)(x + w^2)(x + w^3)(x + w^{20})(x + w^{22})\\ f_2=x(x + w)(x + w^2)(x + w^7)(x + w^{22})(x + w^{25})}&\tabll{y^4 + (2x^6 + w^6x^5 + w^{23}x^4 + 2x^3 + w^6x^2 + w^{25}x)y^2 + w^4x^{10} + w^9x^9 + w^{22}x^8 + w^{10}x^7 + x^6 + w^{24}x^5 + w^{10}x^4 + w^{20}x^3 + w^{20}x^2}&$ 68 $&$ 68-76 $\\
\cline{2-6}
&$6 $&\tabl{f_1=x(x+w)(x+w^2)(x+w^4)(x+w^{20})(x+w^{22})\\ f_2=x(x+w)(x+w^2)(x+w^5)(x+w^{12})(x+w^{24})}&\tabll{y^4 + (2x^6 + x^5 + w^4x^4 + w^{18}x^3 + w^{24}x^2 + w^{10}x)y^2 + w^2x^{10} + w^7*x^9 + w^8x^8 + w^{12}x^7 + w^4x^6 + w^{16}x^5 + w^5x^4 + w^{22}x^3 + w^{24}x^2}&$ 76 $&$ 76-86 $\\
\cline{2-6}
&$7 $&\tabl{f_1=x(x^3+w^{19}x^2+w^{18}x+w^{16})(x^2+w^2x+w^{22})\\ f_2=(x^4+w^7x^2+x+w^4)(x^2+w^2x+w^{22})}&\tabll{y^4 + (2x^6 + w^7x^5 + w^{22}x^4 + w^6x^3 + wx^2 + w^8x + 1)y^2 + w^{12}x^{10} + w^8x^9 + w^{15}x^8 + x^7 + 2x^6 + w^9x^5 + 2x^3 + w^{17}x^2 + 2x + 1}&$ \textbf{82} $&$ 79-96$\\
\cline{2-6}
&$8 $&\tabl{f_1=(x+w^{15})(x^2+w^{25}x+w^5)(x^5+w^3x^4+w^{20}x^3+w^{21}x^2+w^{11}x+w^{18})\\ f_2==(x+w^{18})(x^2+w^2x+w^{22})(x^5+w^3x^4+w^{20}x^3+w^{21}x^2+w^{11}x+w^{18})}&\tabll{y^4 + (2x^8 + wx^7 + w^7x^6 + w^{16}x^5 + w^3x^4 + w^3x^3 + w^{18}x^2 + w^{15}*x + w^{17})y^2 + w^8x^{14} + w^{21}x^{13} + w^{10}x^{12} + w^{23}x^{11} + w^{22}x^{10} + w^{16}x^9 + w^3x^8 + wx^7 + w^{25}x^6 + w^{20}x^5 + w^{15}x^4 + w^{19}x^3 + w^3x^2 + x + w^6}&$ 84 $&$ 92-106$\\
\cline{2-6}
&$9 $&\tabl{f_1=x+w^{12}\\ f_2=x(x^2+w^8x+w^{18})\\ f_3=x^3+w^7x+w^{25}}&\tabll{y^8 + (x^3 + w^{21}x^2 + w^7x)y^6 + (x^6 + w^8x^5 + w^5x^4 + w^{25}x^3 + w^6x^2 + w^4x + w^{11})y^4 + (w^{22}x^7 + w^3x^6 + w^{18}x^5 + w^3x^4 + wx^3 + w^8x^2 + w^{15}*x)y^2 + w^2x^8 + w^{14}x^7 + w^{11}x^6 + w^{24}x^4 + w^8x^3 + w^{23}x^2 + w^2x + w^{22}}&$ 88 $&$ 91-116$\\
\cline{2-6}
&$10$&\tabl{f_1=(x+w)(x^2+w^{12}x+w^{17})(x^5+w^5x^4+x^3+w^{20}x^2+w^8x+w^2) \\ f_2=(x+w)(x^2+w^{12}x+w^{17})(x^5+w^4x^4+wx^3+w^9x^2+w^9x+1) }&\tabll{y^4 + (2x^8 + w^{23}x^7 + w^{14}x^6 + w^2x^5 + w^8x^4 + w^{12}x^3 + w^5x^2 + w^3x + 2)y^2 + w^{14}x^{14} + w^5x^{13} + w^5x^{12} + w^{21}x^{11} + wx^{10} + w^{25}x^9 + w^{15}x^8 + w^{18}x^7 + w^{24}x^6 + w^{22}x^5 + w^{22}x^4 + w^4x^3 + wx^2 + w^{16}x + w1^8}&$ \textbf{94} $&$ 91-126$\\
\cline{2-6}
&$11 $&\tabl{f_1=(x+w)(x^2+w^2)\\ f_2=(x+w^5)(x^2+w^{12}x+w^3)\\ f_3=x^2+w^{18}x+w^{25} }&\tabll{y^8 + (x^3 + w^9x^2 + w^{23}x + w^{24})y^6 + (x^6 + w^{14}x^5 + 2x^4 + w^{24}x^3 + w^3x^2 + 2x + 1)y^4 + (x^8 + w^{11}x^7 + w^{25}x^6 + w^{12}x^5 + w^7x^4 + w^6x^3 + w^6x^2 + w^{17}x + w^{17})y^2 + x^{10} + w^8x^9 + w^{21}x^8 + w^{22}x^7 + w^{14}x^6 + w^{23}x^5 + w^7x^4 + w^{21}x^3 + w^7x^2 + w^{18}x + w^6}&$ \textbf{96} $&$ \leq 136$\\
\hline
\hline
$81$&$1 $&\tabl{f_1=x^2+w, f_2=x^2+w^{41} }&\tabll{y^4 + 2x^2y^2 + w^2}&$ 100 $&$ 100 $\\
\cline{2-6}
&$2 $&\tabl{f_1=(x+w^{57})(x^3+w^2x+w^6)\\ f_2=(x+w^{79})(x^3+w^2x+w^6)}&\tabll{y^4 + (2x^4 + w^{24}x^3 + w^{42}x^2 + w^{16}x + w^{30})y^2 + w^2x^6 + w^{44}x^4 + w^{48}x^3 + w^6x^2 + w^{50}x + w^{14}}&$ 116$ &$ 118 $\\
\cline{2-6}
&$4 $& &\tabll{y^4 + (2x^3 + w^{42}x^2 + w^6x + w^{53})y^2 + w^{48}x^2 + w^{12}x + w^{56}}&$ 150$ &$ 154 $\\
\cline{2-6}
&$5 $&\tabl{f_1=x(x+w^{12})\\ f_2=(x+w^{50})(x+w^{67})\\ f_3=(x+w^4)(x+w^{23}) } &\tabll{y^8 + (w^{50}x + w^7)y^6 + (w^{50}x^3 + w^{22}x^2 + w^{34}x + w^{64})y^4 + (2x^6+ w^{50}x^5 + w^{50}x^4 + w^{27}x^3 + w^8x^2 + w^{26}x + w^{41})y^2 + w^{20}x^6 + w^{55}x^5 + 2x^4 + w^{49}x^3 + w^{25}x^2 + w^4x + w^{68}}&$ \textbf{160}$ &$ 156-172 $\\
\cline{2-6}
&$7 $&\tabl{f_1=(x^2 + w)(x^4 + w^4x^2 + w^{13})\\ f_2=(x^2 + w)(x^4 + w^2x^2 + 2x + w^{22})}&\tabll{y^4 + (2x^6 + w^{12}x^4 + 2x^3 + w^{24}x^2 + w^{41}x + w^{26})y^2 + w^{54}x^8 + w^{67}x^7 + w^3x^6 + w^{51}x^5 + w^{73}x^4 + w^{12}x^3 + w^{50}x^2 + w^{68}x + w^{54}}&$ \textbf{172} $&$ 160-208 $\\
\cline{2-6}
&$13$ &\tabl{f_1=x^4+w^7x^2+w^{77}+w^{16}\\ f_2=x^4+w^{15}x^2+w^{21}x+w^{48} \\ f_3=x^2+w^{23}+w^{59} }&\tabll{y^8 + (x^4 + w^{39}x^2 + w^{15}x + w^{35})y^6 + (x^8 + w^{27}x^6 + w^{9}x^5 + w^{57}x^4 + w^{29}x^3 + w^{65}x^2 + w^{10}x + w^{65})y^4 + (x^{10} + w^{23}x^9 + w^{61}x^8 + w^{12}x^7 + w^{11}x^6 + w^{6}x^5 + w^{22}x^4 + w^{30}x^3 + w^{75}x^2 + w^{38}x + w^{70})y^2 + x^{12} + w^{63}x^{11} + w^{52}x^{10} + w^{46}x^9 + w^{43}x^8 + w^{23}x^7 + w^{68}x^6 + w^{50}x^5 + w^{50}x^4 + w^{75}x^3 + w^{16}x^2 + w^{8}} &$ \textbf{224} $&$ \leq 315$\\
\hline
\end{longtable}
\end{landscape}
\bibliography{liste,zhang,lebrigand,macrae,angles,rosen,aubry,hayes,clement,xing,lachaud,perret,serre,geer,stichtenoth}
\bibliographystyle{plain}
\end{document}